\documentclass[10pt]{amsart}
\usepackage[all]{xy}
\usepackage{amsmath}
\usepackage{amssymb}
\usepackage{amsthm}
\usepackage{parskip}
\input hart.sty         
\renewcommand{\>}{\rightarrow}

\begin{document}
\frenchspacing
\theoremstyle{plain}
\newtheorem{thm}{Theorem}
\newtheorem{prop}[thm]{Proposition}
\newtheorem{cor}[thm]{Corollary}
\newtheorem{lem}[thm]{Lemma}

\theoremstyle{definition}
\newtheorem{defn}[thm]{Definition}
\newtheorem{que}[thm]{Question}
\newtheorem{cla}[thm]{Claim}
\newtheorem{exam}[thm]{Example}
\newtheorem{obs}[thm]{Observation}
\newtheorem{rmk}[thm]{Remark}

\theoremstyle{remark}
\newtheorem{note}[thm]{Note}

\title{Second Infintesimal Neighborhoods of Projective Bundle Sections} 
\author{Travis Kopp}  
\email{travis.kopp@gmail.com}  
\date{\today}       
\maketitle
In this paper we will show that given a smooth variety, $B$, 
there is a correspondence between projective bundles over $B$ with a section 
and second infinitesimal neighborhoods of $B$ of the proper dimension.
(The terminology and background of this paper will be as is found in \cite{Hart1}.)
In particular, we will show that a projective bundle over $B$ with a section can be uniquely recovered
from the embedding of $B$ in its second infinitesimal neighborhood as a section of the bundle.
This follows easily from the following main theorem.

\begin{thm}[Main]\label{main}
Let $B$ be a smooth variety and let the surjection of locally free sheaves,
$\sE\epito\sO_B$, 
determine a projective bundle section, $\si:B\>\LA=\P(\sE)$,
where $\sE$ is a rank $d+1$ locally free sheaf.
Let $\sI_B$ be the ideal sheaf of $\si(B)$ in $\LA$.
Then the following short exact sequences of locally free sheaves on $B$ are isomorphic,
\[
0\>\sK\>\sE\>\sO_B\> 0 \]
\[
0\>\sI_B / \sI_B^2\>\sO_{\LA}/\sI_B^2\>\sO_B\>0 \]
\end{thm}
We will start by describing the situation in which we are interested.
Let $B$ be a smooth variety. Let $\sE$ be a rank $d+1$ locally free sheaf on $B$
with corresponding projective bundle, $\pi:\LA=\P(\sE)\>B$.
A section of this bundle, $\si:B\>\LA$ is determined by a surjection,
$\sE\epito\sL$, of $\sE$ onto an invertible sheaf $\sL$.
By twisting by $\sL^{-1}$, we may assume that $\sL=\sO_B$ (since $\P(\sE\otimes\sL^{-1})\iso\P(\sE)$). 
Let $\sK$ be the kernel of the above surjection of sheaves.
Then the section $\si:B\>\LA$ will correspond 
to the following short exact sequence of locally free sheaves on $B$.
\[
0\>\sK\>\sE\>\O_B\>0 \]
\indent
Let $\sI_B\sub\sO_{\LA}$ be the ideal sheaf of $\si(B)\sub\LA$
and let $\sO_{\LA}(1)$ be the usual relatively very ample invertible sheaf on $\LA$.
We start with the following proposition.

\begin{prop}
\[
\pi_*(\sI_B^2(1))=R^1\pi_*(\sI_B^2(1))=0\]
\end{prop} 

\begin{proof}
We will work locally on $B$. Let $\Spec A \iso U\sub B$
be an open affine neighborhood on which $\sE$ is free.
Then we will have the following correspondences.
\begin{enumerate}
\item 
$\sE\vert_U$ with a free module, $A[x_0,\ldots,x_d]$,
\item
$\pi^{-1}(U)\sub\LA$ with $\Proj A[x_0,\ldots,x_d]=\P_A^d$,
\item
$\sE\vert_U\epito\sO_U$ with a surjection, $A[x_0,\ldots,x_d] \epito A,\ x_i\mapsto s_i$,
\item
$\sI_B\vert_{\pi^{-1}(U)}$ with $I_B=(s_jx_i-s_ix_j)\sub A[x_0,\ldots,x_d]$.
\end{enumerate}

\indent 
We first wish to describe the sheaf $\sO_{\P_A^d}/\sI_B^2(1)$ on the variety $\si(B)$.
Since $A[x_0,\ldots,x_d] \> A,\ x_i\>s_i$ is a surjection,
we can find $\{a_i\}\sub A$ such that $\sum a_is_i = 1$.
Define,
\[
f:=\sum a_ix_i\in A[x_0,\ldots,x_d] \]
\indent
Then $\si(B)\sub D_+(f)\iso \Spec A[x_0,\ldots,x_d]_{(f)}$.
Thus we can work in this open neighborhood when considering $\sO_{\P_A^d}/\sI_B^2(1)$.
In this neighborhood, $\sI_B$ corresponds to,
\[
(I_B)_{(f)}=(s_j\frac{x_i}{f}-s_i\frac{x_j}{f})\sub A[x_0,\ldots,x_d]_{(f)}
\]
Since $\sum a_is_i=\sum a_i\frac{x_i}{f} = 1$, we observe that for a fixed $i$,
\[
\sum_{j}a_j(s_j\frac{x_i}{f}-s_i\frac{x_j}{f})=\frac{x_i}{f}-s_i \]
Therefore $(I_B)_{(f)}$ can be rewritten as,
\[
(I_B)_{(f)}=(y_0,\ldots,y_d) \]
where $y_i=\frac{x_i}{f}-s_i$, with the relation $\sum a_iy_i=0$.\\
\indent
It is clear that $A[x_0,\ldots,x_d]_{(f)}=A[y_0,\ldots,y_d]$.
Therefore $\sO_{\P_A^d}/\sI_B^2$ corresponds to the module,
\[
A[x_0,\ldots,x_d]_{(f)}/(I_B)_{(f)}^2 = A[y_0,\ldots,y_d]/(y_0,\ldots,y_d)^2 \]
and $\sO_{\P_A^d}/\sI_B^2(1)$ corresponds to the module,
\[
\overline{f}(A[y_0,\ldots,y_d]/(y_0,\ldots,y_d)^2) \]
\indent
Consider the following long exact sequence of cohomologies on $\P_A^d$,
\[
0\>H^0(\sI_B^2(1))\>H^0(\sO(1))\> H^0(\sO/\sI_B^2(1))\>H^1(\sI_B^2(1))\>H^1(\sO(1))=0 \]
We know that $H^0(\sO_{\P_A^d}(1))=A[x_0,\ldots,x_d]_1=fA[y_0,\ldots,y_d]_{\leq 1}$.
Therefore the map,
\[
H^0(\sO_{\P_A^d}(1))\>H^0(\sO_{\P_A^d}/\sI_B^2(1)) \]
is an isomorphism of $A$-modules.\\
\indent
From the above exact sequence it then follows that $H^0(\sI_B^2(1))=0$
(which we already knew) and that $H^1(\sI_B^2(1))=0$.
Since this is true locally over $U\sub B$,
the global implication for $\pi:\LA\>B$ is that,
\[
\pi_*(\sI_B^2(1))=R^1\pi_*(\sI_B^2(1))=0\]
\end{proof}
(The following theorem is then the same as Theorem $\ref{main}$.)
\begin{thm}
The following short exact sequences of locally free sheaves on $B$ are isomorphic,
\[
0\>\sK\>\sE\>\sO_B\> 0 \]
\[
0\>\sI_B / \sI_B^2\>\sO_{\LA}/\sI_B^2\>\sO_B\>0 \]
\end{thm}
\begin{proof}
It is clear after a little thought that the short exact sequence,
\[
0\>\sK\>\sE\>\sO_B\> 0 \]
is obtained by applying $\pi_*$ to the natural sequence,
\[
0\>\sI_B(1)\>\sO_\LA(1)\>\sO_B\> 0 \]
(Notice that since $\si$ is determined by $\sE\epito\sO_B$,
we have $\si^*\sO_\LA(1)=\sO_B$.)\\
\indent
Now we consider the diagram,
\[
\begin{array}{ccccccccc}
&&0&&0&&&& \\
&&\downarrow&&\downarrow&&&& \\
0&\>&\sI_B^2(1)&\>&\sI_B^2(1)&\>&0&& \\
&&\downarrow&&\downarrow&&\downarrow&& \\
0&\>&\sI_B(1)&\>&\sO_\LA(1)&\>&\sO_B&\>&0 \\
&&\downarrow&&\downarrow&&\downarrow&& \\
0&\>&\sI_B/\sI_B^2&\>&\sO_\LA/\sI_B^2&\>&\sO_B&\>&0 \\
&&\downarrow&&\downarrow&&\downarrow&& \\
&&0&&0&&0&&
\end{array}
\]
\indent
Since the first two rows are exact and the columns are exact, the third row will also be exact.
Then applying $\pi_*$ to the diagram and taking into account the proposition
yields the statement in the theorem.
\end{proof}
\section*{Acknowledgements}
I would like to thank my Ph.D. advisor, S\'{a}ndor Kov\'{a}cs, 
for his help in proofreading this text and in giving some revision suggestions.

\bibliographystyle{plain}
\bibliography{KodBib}

\end{document}